\documentclass{amsart}

\usepackage{amssymb}

\numberwithin{equation}{subsection}

\newtheorem{theorem}{Theorem}[section]

\theoremstyle{definition}

\theoremstyle{remark}

\newcommand{\apt}{\mathfrak{a}}
\newcommand{\la}{\mathfrak{g}}
\newcommand{\bv}{\mathbf V}
\newcommand{\bla}{\boldsymbol{\mathfrak g}}
\newcommand{\bg}{\mathbf G}
\newcommand{\bF}{\mathbf F}
\newcommand{\bE}{\mathbf E}
\newcommand{\Cartan}{\mathbf t}

\DeclareMathOperator{\im}{im}
\DeclareMathOperator{\Gal}{Gal}
\DeclareMathOperator{\Aut}{Aut}
\DeclareMathOperator{\Ad}{Ad}
\DeclareMathOperator{\Lie}{Lie}
\DeclareMathOperator{\val}{val}
\DeclareMathOperator{\rs}{rs}

\begin{document}

\title{Purity of equivalued affine Springer fibers}
\author{Mark Goresky}\thanks{The research of M.G. was 
supported in part by N.S.F. grants
DMS-9900324 and DMS-0139986}
\author {Robert Kottwitz}\thanks{The research of R.K. was  
supported in part by N.S.F. grant
 DMS-0071971}
\author {Robert MacPherson}
\maketitle

Let $k$ be an algebraically closed field, $G$ a connected reductive group
over~$k$, and $A$ a maximal torus in~$G$. We write $\mathfrak g$ for the Lie
algebra of~$G$ and $\mathfrak a$ for $X_*(A) \otimes_\mathbf Z \mathbf R$. 

Let $F=k((\epsilon))$ be the field of formal Laurent series over~$k$ and let
$\mathfrak o=k[[\epsilon]]$ be the subring of formal power series. We fix an
algebraic closure $\bar F$ of~$F$. We write
$\bg$ and $\bla$ for $G(F)$ and $\mathfrak g(F):=\mathfrak g\otimes_k F$
respectively. Let $y \in \mathfrak a$ and write $\bg_y$ and $\bla_y$ for the
associated (connected) parahoric subgroup and subalgebra respectively. For
example, if $y =0$, then $\bg_y=G(\mathfrak o)$ and $\bla_y=\mathfrak
g(\mathfrak o)$, while if $y$ lies in the interior of an alcove, then $\bg_y$
(resp., $\bla_y$) is an Iwahori subgroup (resp. subalgebra). 

We write $\mathcal F_y$ for the $k$-ind-scheme $\bg/\bg_y$. When $y=0$ for
example, $\mathcal F_y$ is the affine Grassmannian $G(F)/G(\mathfrak o)$. 
Any
$u
\in
\bla$ determines a closed subset
\begin{equation*}
\mathcal F_y(u):=\{g \in \bg/\bg_y : \Ad(g^{-1})(u) \in \bla_y \}
\end{equation*}
of $\mathcal F_y$, 
called an affine Springer fiber. Affine Springer fibers were first studied by
Kazhdan and Lusztig \cite{kazhdan-lusztig88}, who established some of their
basic properties, one of which we will now recall.

Let $T$ be a maximal $F$-torus in~$G$, let $\mathfrak t$ be its Lie algebra
(over
$F$), and let $A_T$ denote the maximal $F$-split torus in~$T$. A regular
element $u \in \mathfrak t$ is said to be \emph{integral} if  $\val \lambda
(u) \ge 0$ for every  $\lambda \in X^*(T)$ (convention: $\val 0=+\infty$).
Here we have denoted the value of the differential of~$\lambda$ on~$u$ simply
by
$\lambda(u)$, an especially abusive notation when $k$ has characteristic $p$,
since in that case the differential of
$\lambda$ is $0$  whenever 
$\lambda$ is divisible by $p$ in $X^*(T)$.  The centralizer
$T(F)$ of~$u$ acts on $\mathcal F_y(u)$, and inside the centralizer we have
the lattice
$\Lambda:=X_*(A_T)$, viewed as a subgroup of~$T(F)$ by sending $\mu \in
\Lambda$ to $\epsilon^\mu \in T(F)$. The lattice $\Lambda$ acts freely on
$\mathcal F_y(u)$, and Kazhdan and Lusztig showed that the quotient
$\Lambda\backslash \mathcal F_y(u)$ is a projective (usually singular)
variety over $k$, non-empty if and only if $u$ is integral, and they
conjectured a formula for its dimension, later proved by Bezrukavnikov
\cite{bezrukavnikov96}.

In order to state the main result of this paper, we need two more 
definitions. We normalize the valuation on $\bar F$ so that $\val \epsilon
=1$.  We say that regular
$u
\in
\mathfrak t$ is
\emph{equivalued with valuation $s \in \mathbf Q$} if  $\val \alpha (u)=s$
for every root $\alpha$ of $T$ over
$\bar F$ and $\val \lambda (u) \ge s$ for every $\lambda \in X^*(T)$. (For
adjoint groups $G$ the second part of the condition is of course redundant.) 
When
$s$ plays no role, we say simply that $u$ is \emph{equivalued}.

Hessenberg varieties were introduced by De Mari, Procesi and Shayman
\cite{demari-procesi-shayman92}. In 
\S~ 1 below we define them in  greater generality but still call them
Hessenberg varieties. The ones in
\cite{demari-procesi-shayman92} are non-empty and admit pavings by affine
spaces, but this is not always true of the more general ones considered
below (see \ref{hempty} for the emptiness question). In any case Hessenberg
varieties are certain very special non-singular closed subvarieties of
(partial) flag manifolds. By a
\emph{Hessenberg paving} of an ind-scheme
$X$ we mean a exhaustive increasing filtration $X_0 \subset X_1 \subset X_2
\subset
\dots$ of $X$ by closed subschemes $X_i$ such that each successive 
difference $X_i \setminus
X_{i-1}$ is a disjoint union of iterated affine space bundles over 
Hessenberg varieties.

\begin{theorem}
Assume that the order of the absolute Weyl group $W$ of $G$ is invertible
in~$k$. Let $u$ be a regular integral equivalued element of $\mathfrak t$.
Then the affine Springer fiber $\mathcal F_y(u)$ admits a Hessenberg paving. 
\end{theorem} 

For some tori $T$ all the (nonempty) Hessenberg varieties turn out to be
single points, in which case the Hessenberg paving becomes a paving by affine
spaces. Before stating this as a precise theorem, we need to recall that
(assuming $|W|$ is invertible in $k$) the $G(F)$-conjugacy classes of
maximal $F$-tori in~$G$ are parametrized by conjugacy classes in~$W$. (This
well-known fact is reviewed in \ref{sect4.2}.) We say that $T$ is
\emph{Coxeter} if the associated conjugacy class in~$W$ consists of Coxeter
elements.  Now consider the centralizer
$M$ of~$A_T$ in~$G$. Then $M$ is a Levi subgroup of~$G$ (hence is split over
$F$) and contains $T$ as an elliptic maximal $F$-torus.  We say that $T$ is
\emph{weakly Coxeter} in $G$ if $T$ is Coxeter  in~$M$. For
example split maximal tori are weakly Coxeter, and all maximal tori in
$GL(n)$ are weakly Coxeter.

Our Hessenberg pavings are obtained by intersecting the affine Springer fiber
with the orbits of a  parahoric subgroup that depends on~$T$. In
\ref{sub3.8} we give a formula for the dimensions of these intersections. 
When the parahoric subgroup turns out to be an Iwahori subgroup, for example
when $T$ is weakly Coxeter, we get a paving by affine spaces, reminiscent of
the one found by De~Concini, Lusztig and Procesi
\cite{deconcini-lusztig-procesi88} for ordinary (that is, non-affine) Springer
fibers, and coinciding with the one found by Lusztig and Smelt
\cite{lusztig-smelt91} in the case of homogeneous $u \in \mathfrak t$, with
$T$ a Coxeter torus in
$GL(n)$.  More precisely, we have the following theorem, in which the
dimensions of the affine spaces mentioned in its statement can be determined
using 
\ref{sub3.8}. 

\begin{theorem}
Assume that  $|W|$  is invertible
in~$k$ and that $T$ is weakly Coxeter. Let $u$ be a regular integral
equivalued element of
$\mathfrak t$. Then the affine Springer fiber $\mathcal F_y(u)$ admits a
 paving by affine spaces. 
\end{theorem}

It is worth noting that the equivalued hypothesis used here is weaker than the
customary one of homogeneity. However, the homogeneity hypothesis
has the virtue of guaranteeing the existence of a useful $\mathbf G_m$-action
on the affine Springer fiber, stemming from the automorphisms $\epsilon
\mapsto z\epsilon$ of~$F$  ($z \in \mathbf G_m$). The work of
Lusztig and Smelt has been generalized to homogeneous elements in $\mathfrak
t$ for Coxeter tori $T$ by Fan \cite{fan96}, Sage \cite{sage00} and Sommers
\cite{sommers97}.  Quite recently these results have been
reconsidered in the context of Cherednik's double-affine Hecke algebras by
Vasserot
\cite{vasserot.pre1} and by Berest, Etingof and Ginzburg
\cite{berest-etingof-ginzburg.pre1}. A homogeneous (weakly Coxeter) case 
for $GL(n)$ was analyzed by
Laumon and Waldspurger \cite{lw}.

It is also worth noting that  Hessenberg pavings which are not pavings by
affines are plentiful for groups other than $GL(n)$. Indeed this occurs in
the example of Bernstein and Kazhdan, discussed in the appendix to
\cite{kazhdan-lusztig88}, in which one of the irreducible components of the
affine Springer fiber for a certain kind of homogeneous  element in 
$\mathfrak{sp}(6)$ admits a dominant map to an elliptic curve.

In the body of the paper we actually allow a more general set-up,
with the role of $ u \in \bla$ being played by 
vectors $v \in \bv:=V \otimes_k F$, where $V$ is a finite dimensional
representation of~$G$ over $k$, and with the role of $\bla_y$ being played by
lattices $\bv_{y,t}$ in~$\bv$ ($t \in \mathbf R$) analogous to the Moy-Prasad
\cite{moy-prasad94} lattices $\bla_{y,t}$ in~$\bla$. In this situation there
are (generalized) affine Springer fibers $\mathcal F_y(t,v)$ (see
\ref{def.asf}), having Hessenberg pavings under suitable hypotheses on $t$
and $v$ (see
\ref{sub3.2}). 

Thus, even in the traditional case $V=\mathfrak g$, in the body of the paper
the parahoric subalgebra $\bla_y$ occurring in the definition of $\mathcal
F_y(u)$ is replaced by a Moy-Prasad  lattice $\bla_{y,t}$,
and we construct a Hessenberg paving of the (generalized) affine Springer
fiber
$\mathcal F_y(t,u)$ for any equivalued $u \in \mathfrak t$ with valuation $s
\ge t$, again under the assumption that $|W|$ is invertible in~$k$.  The
special case
$t=0$ yields the theorem stated earlier.

The main theorem of this paper yields purity results used
to study orbital integrals in 
\cite{gkm.pre1} and in a recent preprint of Laumon \cite{laumon.pre1}. Indeed,
suppose that
$G$ is defined and split over a subfield $k_0$ of~$k$, and that $T$ and $u$
are defined over $k_0((\epsilon))$. Then our Hessenberg pavings are defined
over some finite extension of $k_0$ in~$k$. Moreover it is clear that
when an ind-scheme $S$ over a finite field $\mathbf F_q$ admits a Hessenberg
paving defined over~$\mathbf F_q$, then all eigenvalues of the geometric 
Frobenius on the Borel-Moore homology $H_i^{BM}(X,\mathbf Q_l)$ (with
$l$ prime to $q$) have complex absolute value $q^{i/2}$.

\section{Hessenberg varieties}\label{sec1} Throughout this section we
consider finite dimensional representations $V$ of $G$. For $\lambda \in
X^*(A)$ we let $V_\lambda$ denote the weight space $\{ v \in V :
av=\lambda(a)v \text{ for all $a
\in A$}\}$. 
\subsection{Filtrations on representations $V$ of $G$}\label{sub1.1} Let $y
\in \apt$. As usual $y$ gives rise to an $\mathbf R$-grading on  $V$:
$$ V=\bigoplus_{t \in \mathbf R} V(y,t),
$$
 where $V(y,t)$ is the direct sum of the weight spaces $V_\lambda$ for
weights $\lambda
\in X^*(A)$ such that $\lambda(y)=t$. We will also need the associated
filtrations $F^t_y$ on
$V$, where the subspace $F^t_yV$ is by definition $\oplus_{t' \ge t}V(y,t')$.
Clearly $F^t_yV
\subset F^s_yV$ if $s \le t$. 

\subsection{Unstable vectors} \label{unstab.vect}
 As usual in geometric invariant theory we say
that a vector $v \in V$ is \emph{$G$-unstable} if there exist $g \in G$, $y
\in \apt$ and $t>0$ such that $v \in g\cdot F^t_yV$. 

\subsection{Good vectors} \label{good.vect}
We say that a vector $v \in V$ is \emph{$G$-good} if
there is no non-zero $G$-unstable vector $v^* \in V^*$ that vanishes
identically on the subspace $\la \cdot v$ of $V$. (Here $V^*$ denotes the
vector space dual to $V$.) For example, any regular semisimple element in
$\la$ is $G$-good, at least if the characteristic of~$k$ is different
from~$2$ (see \ref{good.proof}).

\subsection{Partial flag manifold $\mathcal P_y$}Let $y \in \apt$. As usual
$y$ determines a parabolic subgroup $P_y$ of $G$ whose Lie algebra is
$F^0_y\la$. We write $\mathcal P_y$ for the homogeneous space $G/P_y$. The
subspaces $F^t_yV$ are $P_y$-stable, so there is a $G$-equivariant vector
bundle $E^t_yV$ over $\mathcal P_y$ whose fiber at $gP_y \in \mathcal P_y$ is
equal to the subspace $gF^t_yV$ of $V$. 
  
\subsection{Hessenberg varieties} \label{sub1.5} The vector bundle $E^t_yV$
is a subbundle of the constant vector bundle $\tilde V$ over $\mathcal P_y$
with fiber $V$, so we may consider the quotient vector bundle
$\tilde V/E^t_yV$ on $\mathcal P_y$. Any vector $v \in V$ determines a global
section of
$\tilde V/E^t_yV$, and we define $\mathcal P_y(t,v)$ to be the zero-set of
this global section, a closed subscheme of $\mathcal P_y$.  Thus
$\mathcal P_y(t,v)=\{g \in G/P_y : g^{-1}v \in F^t_yV \}$. 

The map $Y\mapsto Y\cdot v$ from $\la$ to $V$ can be viewed as a constant
map  $\cdot v:\tilde{\la}
\to \tilde V$ of constant vector bundles. Over $\mathcal P_y(t,v)$ the vector
$v$ lies in $E^t_yV$, so the map $\cdot v$ carries
$E^0_y\la$ into $E^t_yV$, yielding an induced map
$$
\cdot v:\tilde{\la}/E^0_y\la \to \tilde V/E^t_yV.
$$ The tangent bundle to $\mathcal P_y$ is $\tilde{\la}/E^0_y\la$, and it is
easy to see that the Zariski tangent space to $\mathcal P_y(t,v)$ at a point
of $ \mathcal P_y(t,v)$ is equal to the kernel of
$
\cdot v:\tilde{\la}/E^0_y\la \to \tilde V/E^t_yV
$ on the fibers at that point.

Now suppose that $v$ is a $G$-good vector and that $t\le 0$, in which case we
refer to $\mathcal P_y(t,v)$ as a \emph{Hessenberg variety}. We claim that
any Hessenberg variety is non-singular and projective. Projectivity is clear.
To prove non-singularity we will check that the global section
$v$ of $\tilde V/E^t_yV$ is transverse to the zero-section, or in other words
that the map 
$
\cdot v:\tilde{\la}/E^0_y\la \to \tilde V/E^t_yV
$ is surjective on fibers at each point $gP_y \in \mathcal P_y(t,v)$. To
prove surjectivity at this point we must show that any $v^* \in V^*$ that
annihilates both $\la \cdot v$ and $g \cdot F^t_yV$ is zero. Since $v^*$
annihilates the second subspace and $t \le 0$, there exists $\delta > 0$ such
that $v^* \in g\cdot F^\delta_yV^*$, and therefore $v^*$ is $G$-unstable.
Since $v$ is $G$-good, it follows that $v^*=0$.

\subsection{When are Hessenberg varieties empty?}\label{hempty} 
Hessenberg
varieties are sometimes empty. Fortunately there is a fairly practical way to
determine when this happens, as we will now see. With notation as above,
choose a Borel subgroup $B$ of~$G$ containing~$A$ and contained in~$P_y$. The
closed subvariety 
\begin{equation*}
Y:= \{g \in G/B : g^{-1}v \in F^t_yV \}
\end{equation*}
of the flag manifold $G/B$ obviously maps onto the Hessenberg variety
$\mathcal P_y(t,v)$, and therefore it is enough to determine when $Y$ is
empty. 

Now $Y$ is empty precisely when its fundamental class in
$H^\bullet(G/B,\mathbf Q_l)$ is $0$ (where $l$ is a prime that is non-zero in
our ground field $k$). Recall that the $\mathbf Q_l$-algebra
$H^\bullet(G/B,\mathbf Q_l)$ is the quotient of the symmetric algebra
$S^\bullet$ on $X^*(A) \otimes_\mathbf Z \mathbf Q_l$ by the ideal
$I$ generated by the Weyl group invariant elements that are homogeneous of
strictly positive degree. 

It follows from the discussion in \ref{sub1.5} that $Y$ is obtained as the
zero-set of a global section of the $G$-equivariant vector bundle on~$G/B$
obtained from the representation $V/F^t_yV$ of~$B$, and that moreover this
global section is transverse to the zero-section. Therefore the fundamental
class of~$Y$ in $H^\bullet(G/B,\mathbf Q_l)$ is equal to the top Chern class
of our vector bundle on~$G/B$. Filtering $V/F^t_yV$ by a complete flag of
$B$-stable subspaces, we see that our top Chern class is the image under
$S^\bullet \twoheadrightarrow H^\bullet(G/B,\mathbf Q_l)$ of the element
$\lambda_1\dots \lambda_m$, where $m$ is the dimension of $V/F^t_yV$ and
$\lambda_1,\dots,\lambda_m$ are the characters of $A$ (with multiplicities)
occurring in the $B$-module $V/F^t_yV$.

We conclude that our Hessenberg variety is empty if and only if 
$\lambda_1\dots \lambda_m$ lies in the ideal $I$.

\subsection{Some vector bundles over Hessenberg varieties in the graded
case}\label{sub1.6} 
The following constructions will be used  to produce the 
Hessenberg pavings of section 3.   

Let $x
\in
\apt$. We use
$x$ to define an
$\mathbf R/\mathbf Z$-grading  
$$ V=\bigoplus_{s \in \mathbf R/\mathbf Z} V(x,s+\mathbf Z)
$$ on $V$, where 
$$ V(x,s+\mathbf Z):=\bigoplus_{m \in \mathbf Z}V(x,s-m),
$$ with $V(x,s-m)$ as in \ref{sub1.1}. In this way $\la$ itself becomes an
$\mathbf R/\mathbf Z$-graded Lie algebra and $V$ becomes an 
$\mathbf R/\mathbf Z$-graded $\la$-module. 

We denote by $H$ the subgroup of $G$ generated by $A$ and the root subgroups
for all roots
$\alpha$ of $A$ in $\la$ such that $\alpha(x) \in \mathbf Z$. Then $H$ is a
connected reductive group whose Lie algebra is $\la(x,0+\mathbf Z)$. The
action of $H$ on $V$ preserves the 
$\mathbf R/\mathbf Z$-grading of $V$. 

Now let $y \in \apt$ and $t \in \mathbf R$. Since $A$ is also a maximal torus
in $H$, we may apply the construction in \ref{sub1.1} to the representation
$V(x,\mathbf s)$ of $H$, where we have written $\mathbf s$ as an abbreviation
for $s + \mathbf Z$. In this way we get a vector subbundle
$E^t_yV(x,\mathbf s)$ of the constant vector bundle $\tilde V(x,\mathbf s)$ on
$\mathcal P_y=H/P_y$. 

Let us fix $x,y \in \apt$, $s,t \in \mathbf R$ and $v \in V(x,\mathbf s)$. We
assume further that
$t
\le 0$ and that the image of
$v$ in $V$ is a $G$-good vector (which implies that $v$ is $H$-good in
$V(x,\mathbf s)$). Under these hypotheses we are going to define some vector
bundles over the Hessenberg variety
$\mathcal P_y(t,v) \subset \mathcal P_y=H/P_y$ associated to the vector $v$
in the representation
$V(x,\mathbf s)$ of $H$. 

The vector bundles depend on two additional real numbers $r$ and $t'$. We  
write $\mathbf r$ for $r+\mathbf Z$,  and we put $t'':=t+t'$. We assume that
$t''\le 0$. The action of $\la$ on $V$ induces a linear map
\begin{equation}\label{star0}
\la(x,\mathbf r) \otimes V(x,\mathbf s) \to V(x,\mathbf r + \mathbf s).
\end{equation}  This gives rise to a constant homomorphism
$$
\tilde{\la} (x,\mathbf r) \otimes \tilde V(x,\mathbf s) \to \tilde
V(x,\mathbf r + \mathbf s)
$$ of constant vector bundles which carries $E^{t'}_y \la(x,\mathbf r) \otimes
E^t_yV(x,\mathbf s)$ into $E^{t''}_yV(x,\mathbf r + \mathbf s)$. Over the
Hessenberg variety
$\mathcal P_y(t,v)$ the global section $v$ lies in the subbundle
$E^t_yV(x,\mathbf s)$. Therefore the map $Y \mapsto Y \cdot v$ induces a
homomorphism 
\begin{equation}\label{1.6.2}
\cdot v: \tilde{\la} (x,\mathbf r)/E^{t'}_y \la(x,\mathbf r)  \to \tilde
V(x,\mathbf r + \mathbf s)/E^{t''}_yV(x,\mathbf r + \mathbf s),
\end{equation} and in fact this homomorphism is surjective at every point of
the Hessenberg variety, as one sees using the $G$-goodness of $v$ and the
fact that $t'' \le 0$, just as in the proof of non-singularity of Hessenberg
varieties. It follows that the kernel of the homomorphism, call it
$K(x,y,r,s,t,t',v)$, is a vector bundle. Denoting by $H_v$ the stabilizer of
$v$ in $H$ (which of course acts on our Hessenberg variety), we have an exact
sequence of $H_v$-equivariant vector bundles on the Hessenberg variety
$\mathcal P_y(t,v)$:
$$ 0 \to K(x,y,r,s,t,t',v) \to \tilde{\la} (x,\mathbf r)/E^{t'}_y
\la(x,\mathbf r)  \to \tilde V(x,\mathbf r + \mathbf s)/E^{t''}_yV(x,\mathbf
r + \mathbf s) \to 0.
$$

\section{Moy-Prasad gradings and filtrations}
\subsection{Moy-Prasad filtrations on representations $V$ of
$G$}\label{sub2.1} As in the introduction we write $\bv$ for $V\otimes_k F$
and $\bla$ for $\la\otimes_k F$, where
$F=k((\epsilon))$. The subspace $V\otimes_k k[\epsilon,\epsilon^{-1}]$ of
$\bv$ can be written as a direct sum of the subspaces $V_\lambda \epsilon^m$,
where $(\lambda,m)$ ranges through
$X^*(A)\oplus \mathbf Z$. 

Now suppose that $ x \in \apt$. For $r \in \mathbf R$ we define a subspace
$\bv_x(r)$ of
$\bv$ by
$$
\bv_x(r)=\bigoplus_{\lambda(x)+m=r}V_\lambda \epsilon^m.
$$ (In other words the direct sum is taken over all $(\lambda,m) \in
X^*(A)\oplus \mathbf Z$ such that 
$\lambda(x)+m=r$.) Thus we have
$$V\otimes_k k[\epsilon,\epsilon^{-1}]=\bigoplus_{r \in \mathbf R} \bv_x(r).$$
Now put $\bv_{x,r}:=\prod_{r' \ge r} \bv_x(r')$. Clearly
$\bv_{x,r} \subset \bv_{x,s}$ if $r \ge s$, so we have defined an $\mathbf
R$-filtration on
$\bv$.   Put $\bv_{x,r+}:=\prod_{r' > r} \bv_x(r')$. Thus we have exact
sequences
$$ 0 \to \bv_{x,r+} \to \bv_{x,r} \to \bv_x(r) \to 0.
$$ Clearly $\bv_{x,r+}=\bv_{x,r+\delta}$ for all sufficiently small $\delta >
0$.

Taking $V=\la$, we get  $\bla_x(r)$ and $\bla_{x,r}$. Note that $\bla_{x,r}$
are the lattices in $\bla$ defined by Moy-Prasad \cite{moy-prasad94}. In this
way
$\bla$ becomes an $\mathbf R$-filtered Lie algebra over $k$ and $\bv$ an
$\mathbf R$-filtered $\bla$-module. In other words the action of
$\bla$ on $\bv$ induces linear maps
\[
\bla_{x,r} \otimes \bv_{x,s} \to \bv_{x,r+s}
\] and
\begin{equation} \label{star1}
\bla_x(r) \otimes \bv_x(s) \to \bv_x(r+s). 
\end{equation}

Recall the subspace $V(x,\mathbf r)$ of $V$ defined in \ref{sub1.6}, where
$\mathbf r$ denotes $r +
\mathbf Z$, as before.
 We are now going to define a canonical linear isomorphism  
\begin{equation} 
i_{x,r}:\bv_x(r)
\to  V(x,\mathbf r).
\end{equation} 
The space $ \bv_x(r)$ is the direct sum of the spaces
$V_\lambda \epsilon^m$ for $\lambda,m$ such that $\lambda(x)+m=r$. For such
$\lambda,m$ define a linear isomorphism  $V_\lambda \epsilon^m \to V_\lambda$
by sending $v\epsilon^m$ to $v$. The direct sum of these isomorphisms is the
desired isomorphism 
$  \bv_x(r) \to  V(x,\mathbf r)$. Under these isomorphisms the linear map
\eqref{star1} goes over to the linear map \eqref{star0}.

\subsection{Moy-Prasad subgroups of $\bg$}\label{sub2.2} We denote $G(F)$ by
$\bg$. For $x \in \apt$ and $r \ge 0$  Moy-Prasad \cite{moy-prasad94} have 
defined a subgroup $\bg_{x,r}$ of $\bg$. We abbreviate $\bg_{x,0}$ to
$\bg_x$, a parahoric subgroup of $\bg$. Viewing $\bg_{x,r}$ as a
pro-algebraic group over $k$,
 its Lie algebra can be identified with $\bla_{x,r}$. The subgroups
$\bg_{x,r}$ are normal in
$\bg_x$. For $s \ge r$ we have $\bg_{x,s} \subset \bg_{x,r}$. We write
$\bg_{x,r+}$ for $\cup_{s >r}\bg_{x,s}$, a subgroup of $\bg_{x,r}$; for
sufficiently small $\delta > 0$ $\bg_{x,r+\delta}$ is equal to  $\bg_{x,r+}$.
Put $\bg_x(r):=\bg_{x,r}/\bg_{x,r+}$.  

Let $V$ be a representation of $G$ and
let $s$ be a real number.  Then the parahoric subgroup $\bg_x$ preserves
$\bv_{x,s}$. Moreover, for every $r 
 \ge 0$ the subgroup $\bg_{x,r}$ acts trivially on $\bv_{x,s}/\bv_{x,s+r}$.
Finally, for $r > 0$ there is a canonical isomorphism 
\[
\xi_{x,r}:\bg_x(r) \to \bla_x(r),
\] and this isomorphism has the following property. Let $g \in \bg_{x,r}$ and
let $v \in
\bv_{x,s}$. Then $gv-v$ lies in $\bv_{x,r+s}$, and its image $w$ in
$\bv_x(r+s)$ depends only on the image $\bar g$ of $g$ in $\bg_x(r)$ and the
image $v_s$ of $v$ in $\bv_x(s)$. The property of $\xi_{x,r}$ we are
referring to is that
\begin{equation} \label{star2} w=\xi_{x,r}(\bar g)\cdot v_s
\end{equation} (The right-hand side of this equation uses the map
\eqref{star1} induced by the action of $\bla$ on
$\bv$.)

\subsection{More filtrations and vector bundles}\label{sub2.3} We will now
see that the affine set-up of \ref{sub2.1}, \ref{sub2.2} produces Hessenberg
varieties and vector bundles of the type studied in  \ref{sub1.6}.

Suppose we have $x,y \in \apt$ and $r,t \in \mathbf R$. We define a subspace
$\bF^t_y\bv_x(r)$ of
$\bv_x(r)$ by  
\[
\bF^t_y\bv_x(r):=\im[\bv_{y,t} \cap \bv_{x,r} \to
\bv_{x,r}/\bv_{x,r+}=\bv_x(r)].
\] Applying this construction to $V=\la$ and $t=0$, we get
\[
\bF^0_y\bla_x(r):=\im[\bla_{y,0} \cap \bla_{x,r} \to
\bla_{x,r}/\bla_{x,r+}=\bla_x(r)].
\] Analogously, for $r \ge 0$ we have a subgroup
\[
\bF^0_y\bg_x(r):=\im[\bg_{y} \cap \bg_{x,r} \to
\bg_{x,r}/\bg_{x,r+}=\bg_x(r)],
\] and for $r > 0$ the subgroup $\bF^0_y\bg_x(r)$ goes over to
$\bF^0_y\bla_x(r)$ under the isomorphism $\xi_{x,r}:\bg_x(r) \to \bla_x(r)$.

Put $H_0:=\bg_x(0)$, a connected reductive group over $k$, and put
$P=\bF^0_y\bg_x(0)$, a parabolic subgroup of $H_0$. Each space $\bv_x(r)$ is
a representation of $H_0$, and $P$ preserves the subspace $\bF^t_y\bv_x(r)$.
Therefore, as in section \ref{sec1}, we may define an
$H_0$-equivariant vector bundle $\bE^t_y\bv_x(r)$ over $\mathcal P:=H_0/P$
whose fiber over $hP$ is equal to $h\bF^t_y\bv_x(r)$.

To see that these constructions are instances of ones made in 
\ref{sub1.6}, let us examine $\bF^t_y\bv_x(r)$ more closely. This subspace of
$\bv_x(r)$ is the direct sum of the spaces $V_\lambda\epsilon^m$ for
$(\lambda,m)$ such that
$\lambda(x)+m=r$ and $\lambda(y)+m\ge t$. These two conditions can be
rewritten as 
$\lambda(x)+m=r$ and $\lambda(y-x)\ge t-r$, which shows that under the
isomorphism $i_{x,r}$ from $\bv_x(r)$ to $V(x,\mathbf r)$, the subspace
$\bF^t_y\bv_x(r)$ goes over to
$\bF^{t-r}_{y-x}V(x,\mathbf r)$. In particular, under the isomorphism
$\bla_x(r)=\la(x,\mathbf r)$, the subspace $\bF^0_y\bla_x(r)$ goes over to
$F^{-r}_{y-x}\la(x,\mathbf r)$.

Note that $H_0$ can be identified with the group $H$ defined in \ref{sub1.6}
(for the given $x \in
\apt$), and that this identification is compatible with the identification
$\bla_x(0)=\la(x,\mathbf 0)$ of their Lie algebras. Under this
identification, the parabolic subgroup $P \subset H_0$ goes over to $P_{y-x}
\subset H$, and the representation $\bv_x(r)$ of
$H_0$ goes over to the representation $V(x,\mathbf r)$ of $H$. We may also
identify $\mathcal P=H_0/P$ with $\mathcal P_{y-x}$, and then the vector
bundles $\bE^t_y\bv_x(r)$, $\bE^0_y\bla_x(r)$ go over to
$E^{t-r}_{y-x}V(x,\mathbf r)$, $E^{-r}_{y-x}\la(x,\mathbf r)$ respectively. 

\section{Generalized affine Springer fibers} 
\subsection{Definition of generalized affine Springer fibers} \label{def.asf}
 Fix $y \in
\apt$. We write $\mathcal F_y$ for the (partial) affine flag space
$\bg/\bg_y$, an ind-scheme over $k$. For $t \in \mathbf R$ and $v \in \bv$ we
define the generalized affine Springer fiber $\mathcal F_y(t,v) \subset
\mathcal F_y$ by
\[
\mathcal F_y(t,v)=\{g \in \bg/\bg_y:g^{-1}v \in \bv_{y,t}\}, 
\] a closed subset of $\mathcal F_y$. 
 
\subsection{Hypotheses on $t$ and $v$}\label{sub3.2} We are going to analyze
the structure of this generalized affine Springer fiber under the following
hypotheses on $t,v$. We suppose that there exist $x \in
\apt$ and $s \in \mathbf R$ such that
\begin{enumerate}
\item $s \ge t$, 
\item $v \in \bv_{x,s}$, and
\item the image $\bar v$ of $v$ under $\bv_{x,s} \to \bv_x(s) \to V(x,\mathbf
s)
\hookrightarrow V$ is a $G$-good vector in $V$. 
\end{enumerate} In the third condition the map $\bv_{x,s} \to \bv_x(s)$ is
the canonical projection and the map
$\bv_x(s) \to V(x,\mathbf s)$ is the linear isomorphism $i_{x,s}$ defined in
\ref{sub2.1}.  
\subsection{Intersection with orbits} In order to analyze  $\mathcal
F_y(t,v)$ we intersect it with the orbits of $\bg_x$ in $\mathcal F_y$. The
intersections will turn out to be iterated affine space bundles over
Hessenberg varieties (and are therefore empty precisely when the Hessenberg
variety is empty).   

We first make the observation that it is enough to
consider the orbit of the base-point in
$\mathcal F_y$. Indeed, right multiplication by an element $c$ in the
normalizer of $A(F)$ in $G(F)$ induces an isomorphism from the intersection
of $\mathcal F_{cy}(t,v)$ and the $\bg_x$-orbit of the base-point in
$\mathcal F_{cy}$ to the intersection of $\mathcal F_y(t,v)$ and the
$\bg_x$-orbit of $c$ times the base-point in $\mathcal F_y$. (Recall that $c$
acts on $\apt$ by an affine linear transformation; $cy$ denotes the result of
this action on~$y$.)  

The intersection of $\mathcal F_y(t,v)$ and the
$\bg_x$-orbit of the base-point in $\mathcal F_y$ is equal to
\[ S:=\{ g \in \bg_x/(\bg_x \cap \bg_y) : g^{-1}v \in \bv_{y,t}\}.
\]

\subsection{Spaces $\tilde S_r$ and $S_r$} We must analyze $S$, and in order
to do so we introduce some auxiliary spaces. For $r > 0$ put
\[
\tilde S_r:=\{ g \in \bg_x/(\bg_x \cap \bg_y) : g^{-1}v \in
\bv_{y,t}+\bv_{x,s+r} \}
\] and for $r \ge 0$ put
\[
\tilde S_{r+}:=\{ g \in \bg_x/(\bg_x \cap \bg_y) : g^{-1}v \in
\bv_{y,t}+\bv_{x,(s+r)+} \}.
\] Note that $\tilde S_{r+}=\tilde S_{r+\delta}$ for all sufficiently small
$\delta>0$. 

Since $v \in \bv_{x,s}$ we have $g^{-1}v \in \bv_{x,s}$ for any $g \in
\bg_x$, and the condition
$g^{-1}v \in \bv_{y,t}+\bv_{x,s+r}$ can be thought of as the condition that
the image of $g^{-1}v$ in $\bv_{x,s}/\bv_{x,s+r}$ lie in the subspace
\[
\im[\bv_{y,t}\cap \bv_{x,s} \to \bv_{x,s}/\bv_{x,s+r}].
\] Since $\bg_{x,r}$ is normal in $\bg_x$ and acts trivially on
$\bv_{x,s}/\bv_{x,s+r}$, we see that
$\tilde S_r$ (respectively, $\tilde S_{r+}$) is invariant under the left
action of
 $\bg_{x,r}$ (respectively,
$\bg_{x,r+}$), so that we may define quotient spaces
\begin{align*} S_r := &\bg_{x,r} \backslash \tilde S_r \\ S_{r+} :=
&\bg_{x,r+} \backslash \tilde S_{r+}
\end{align*}
 Note that $S_{r+}=S_{r+\delta}$ for all sufficiently small
$\delta > 0$. 

Let $0=r_0 <r_1<r_2<r_3<\dots$ be the discrete set of values of $r\ge 0$ for
which  either $\bg_x(r)$ or
$\bv_x(s+r)$ is non-trivial. We have inclusion maps
\[
\dots \to \tilde S_{r_3} \to \tilde S_{r_2} \to \tilde S_{r_1}
\] and these induce maps
\[
\dots \to  S_{r_3} \to  S_{r_2} \to  S_{r_1}.
\]

Note that $S_{r_j}=S$ for sufficiently large $j$, and that
$S_{r_{i+1}}=S_{r_i+}$ for all
 $i \ge 0$ (in particular, $S_{r_1}=S_{0+}$). Therefore, in order to show
that $S$ is an iterated affine space bundle over a Hessenberg variety it is
enough to show that $S_{0+}$ is a Hessenberg variety and that for every
$r>0$ the map $S_{r+} \to S_r$ is an affine space bundle. In fact we will
prove a more precise statement. For each $r \ge 0$ we are going to construct
a vector bundle $E_r$ over the Hessenberg variety $S_{0+}$. It will turn out
that $E_0$ is the tangent bundle to $S_{0+}$ and that $S_{r+}
\to S_r$ is a torsor under the vector bundle on $S_r$ obtained by pulling
back the vector bundle
$E_r$ by means of $S_r \to S_{0+}$. Thus the dimension of $S$ is given by
$\sum_{r \ge 0}\dim E_r$, the sum having only finitely many non-zero terms.
However it is conceptually clearer to calculate the dimension slightly
differently, as we will do in 
\ref{sub3.8} below.

\subsection{Analysis of $S_{0+}$} We now check that $S_{0+}$ is a Hessenberg
variety. We use the notation of \ref{sub2.3}. It is immediate from the
definitions that 
\[ S_{0+}=\{ g \in H_0/P : g^{-1}v_s \in \bF^t_y\bv_x(s) \}
\] where $v_s$ denotes the image of $v$ under the canonical projection
$\bv_{x,s} \to \bv_x(s)$. It follows from the discussion in \ref{sub2.3} that
$S_{0+}$ equals $\mathcal P_{y-x}(t-s,\bar v)$, where $\bar v$ is the image
of $v_s$ under the isomorphism
$i_{x,s}:\bv_x(s) \to V(x,\mathbf s)$. Since $t-s \le 0$ and $\bar v$ is
$G$-good by hypothesis, it follows that $\mathcal P_{y-x}(t-s,\bar v)$ is a
Hessenberg variety. 
\subsection{Some auxiliary bundles}\label{sub3.6} Fix a real number $r > 0$. 
Before analyzing the fibers of $S_{r+} \to S_r$, we introduce some auxiliary
bundles. Consider the map
\[ p:\bg_{x,r+}\backslash \bg_x \to \bg_{x,r}\backslash \bg_x.
\] Since $\bg_x(r)=\bla_x(r)$, it is clear that $p$ is a principal
$\bla_x(r)$-bundle, or in other words a torsor under the constant vector
bundle $\tilde \bla_x(r)$  over $\bg_{x,r} \backslash \bg_x$ with fiber
$\bla_x(r)$. 

Next consider the map
\[ q:\bg_{x,r+}\backslash \bg_x / (\bg_x \cap \bg_y) \to \bg_{x,r}\backslash
\bg_x / (\bg_x \cap
\bg_y).
\] The target of $q$ maps to 
\[
\bg_{x,0+}\backslash \bg_x / (\bg_x \cap \bg_y)=H_0/P=\mathcal P.
\] On $\mathcal P$ we have the vector bundle $\bE^0_y\bla_x(r)$, and it is
easy to see that $q$ is a torsor under the quotient of $\tilde \bla_x(r)$ by
the pullback of $\bE^0_y\bla_x(r)$.

\subsection{Analysis of  $S_{r+} \to S_r$} We continue with $r$ and $q$ as in
\ref{sub3.6}. Of course $S_{r+}$ is a closed subspace of the source of $q$
and $S_r$ is a closed subspace of the target of $q$, and it is the
restriction of $q$ to
$S_{r+}$ that gives the map $S_{r+} \to S_r$ that we need to understand.

In order to see that $S_{r+} \to S_r$ is a bundle, we begin by making the
following simple calculation. Suppose that $g \in \bg_x$ represents an
element of $S_r$, so that $g^{-1}v \in
\bv_{y,t}+\bv_{x,s+r}$. Then $g$ represents an element of $S_{r+}$ if and
only if the image of $v$ is $0$ in
\[
\frac{g\bv_{y,t}+\bv_{x,s+r}}{g\bv_{y,t}+\bv_{x,(s+r)+}}=\frac{\bv_x(s+r)}{g\bF^t_y\bv_x(s+r)}.
\] 

Suppose that we replace $g$ by $g_rg$ with $g \in \bg_{x,r}$. (In other
words we are considering
 points in the fiber of $q$ through the point represented by $g$.) Note that
\[ g_rg\bF^t_y\bv_x(s+r)=g\bF^t_y\bv_x(s+r).
\] Since $v \in \bv_{x,s}$ we have $g_r^{-1}v \equiv v-\xi_{x,r}(g_r) \cdot
v$ modulo
$\bv_{x,(r+s)+}$. Therefore the element of $\bv_x(s+r)/g\bF^t_y\bv_x(s+r)$
obtained from~$g_rg$ differs from the one obtained from~$g$ by the element
$\xi_{x,r}(g_r)\cdot v_s \in \bv_x(s+r)$, where $v_s$ denotes the image
of~$v$ in~$\bv_x(s)$.   Therefore $g \in S_r$ lies in the image of $S_{r+}$
if and only if the image of $v$ in 
\[\frac{\bv_x(s+r)}{g\bF^t_y\bv_x(s+r)}\] lies in the image of
\begin{equation} \label{dim.filt}
\cdot v_s:\bla_x(r) \to \frac{\bv_x(s+r)}{g\bF^t_y\bv_x(s+r)}.
\end{equation}

Using the linear isomorphisms $\bla_x(r) \simeq \la(x,\mathbf r)$ and
$\bv_x(s+r)  \simeq  V(x,\mathbf s + \mathbf r)$ discussed at the end of
\ref{sub2.1},  we can rewrite this last map as
\begin{equation}
\cdot \bar v:\la(x,\mathbf r) \to V(x,\mathbf s + \mathbf r)/\bar g 
F^{t-s-r}_{y-x}V(x,\mathbf s+\mathbf r),
\end{equation} 
where $\bar g$ denotes the image of $g$ in $\bg_x(0)=H_0$. (Recall from
\ref{sub2.3} that $H_0$ can be identified with the group $H$ defined in
\ref{sub1.6}, and that  
$V(x,\mathbf s+\mathbf r)$ is a representation of~$H$.)  Our hypotheses
guarantee that 
$t-s-r <0$ and that $\bar v$ is a $G$-good vector in~$V$. 
It then follows from
\ref{sub1.6} that the map
$\cdot
\bar v$ above is surjective, and hence that $g \in S_r$ does lie in the image
of $S_{r+}$.  

Therefore $S_{r+} \to S_r$ is surjective, and the analysis
above (together with \ref{sub3.6})  shows that the fiber of
$S_{r+}
\to S_r$ through the point in $S_{r+}$ represented by $g \in \bg_x$ is equal
to the kernel of the map $\cdot \bar v$ above modulo $g\bF^0_y\bla_x(r)=\bar
gF^{-r}_{y-x}\la(x,\mathbf r)$.   This shows that
$S_{r+} \to S_r$ is a torsor under the vector bundle on $S_r$ obtained by
pulling back the vector bundle
\begin{align*} E_r:&=\ker[\cdot v_s: \tilde \bla_x(r) /\bE^0_y\bla_x(r) \to
\tilde \bv_x(s+r)/\bE^t_y\bv_x(s+r)]\\
   &=\ker[\cdot \bar v:\tilde {\la}(x,\mathbf r)/E^{-r}_{y-x}\la(x,\mathbf r)
\to \tilde V(x,\mathbf s+\mathbf r)/E^{t-s-r}_{y-x}V(x,\mathbf s + \mathbf
r)],
\end{align*}
 the vector bundle denoted by $K(x,y-x,r,s,t-s,-r,\bar v)$ in
\ref{sub1.6}. Note that the vector bundle $E_r$ is also defined for $r=0$, in
which case it coincides with the tangent bundle  (see \ref{sub1.5}) of the
Hessenberg variety $S_{0+}=\mathcal P_{y-x}(t-s,\bar v)$.

\subsection{Dimension of $S$}\label{sub3.8}
Now we are going to calculate the dimension of $S$.  The group $\bg_x$ acts
on the vector space $\bv_{x,s}$. Over $\bg_x/(\bg_x \cap
\bg_y)$ we have the vector bundle whose fiber at $g \in \bg_x/(\bg_x \cap
\bg_y)$ is the vector space $\bv_{x,s}/(\bv_{x,s} \cap g\bv_{y,t})$, and the
vector $v \in \bv_{x,s}$ determines a section of this vector bundle, a section
whose zero-set coincides with $S$, as follows immediately from the
definitions. 

In fact this section is transverse to the zero-section.  In other words, the
map 
\begin{equation*}
\cdot v: \bla_x \to \bv_{x,s}/(\bv_{x,s} \cap g\bv_{y,t})
\end{equation*}
is surjective for every $g \in \bg_x$ that represents a point in~$S$, as
follows from  the
surjectivity of the maps
\eqref{dim.filt}. 
Transversality implies that
\begin{equation}\label{dim.s}
\dim(S)=\dim(\bla_x/(\bla_x \cap \bla_y))-\dim(\bv_{x,s}/(\bv_{x,s}\cap
\bv_{y,t}).
\end{equation}
It is clear that $\dim(\bla_x/(\bla_x \cap \bla_y))$
is the number of affine roots that are non-negative on $x$ and strictly
negative on~$y$. In the  important special case when $\bv=\bla$ and
$t=0$, there is a similar expression for the last dimension appearing in
equation \eqref{dim.s}, and we conclude that in this special case the
dimension of $S$ is the cardinality of
\begin{equation}
\{ \text{affine roots }\tilde\alpha : 0 \le \tilde\alpha(x) < s \text{ and }
\tilde\alpha(y) < 0 \}.
\end{equation}

\section{Affine Springer fibers in the equivalued case}
\subsection{Additional hypothesis on $k$} Throughout this section we assume
that the order of the Weyl group $W:=N_G(A)/A$ of $G$ is invertible in~$k$. 

\subsection{$\bg$-conjugacy classes of maximal $F$-tori in $G_F$}
\label{sect4.2} We write
$G_F$ for the $F$-group obtained from $G$ by extension of scalars from $k$ to
$F$, and as before we write $\bg$ for $G(F)$. We need to review the
well-known classification of
$\bg$-conjugacy classes of maximal $F$-tori in $G$. 

Let $T$ be a maximal $F$-torus in $G_F$. Fix an algebraic closure $\bar F$
of~$F$ and let $F_s$ denote the separable closure of 
$F$ in~$\bar F$. Choose $h \in G(F_s)$ such that $T=hAh^{-1}$. For $\sigma
\in
\Gamma:=\mathrm{Gal}(F_s/F)$ put
$x_\sigma=h^{-1}\sigma(h)$. Then $\sigma \mapsto x_\sigma$ is a $1$-cocycle
of $\Gamma$ in
$G(F_s)$, whose class in $H^1(F,N_G(A))$ we denote by $x(T)$. This
construction yields a bijection from the set of $\bg$-conjugacy classes of
maximal $F$-tori in $G_F$ to the set
\[
\ker[H^1(F,N_G(A)) \to H^1(F,G)].
\]

The obvious surjection $N_G(A) \to W$ induces an injection
\[ H^1(F,N_G(A)) \to H^1(F,W)
\] since for our particular field $F$ it is known \cite{serre68} that
$H^1(F,T')$ vanishes for every
$F$-torus $T'$. Denoting by $y(T)$ the image of $x(T)$ in $H^1(F,W)$, we
obtain an injection $T
\mapsto y(T)$ from the set of  $\bg$-conjugacy classes of maximal $F$-tori in
$G_F$ to
$H^1(F,W)$. 

In fact all elements of $H^1(F,W)$ arise in this way, and we will now review
a proof of this fact, since for our purposes we need to produce an explicit
maximal
$F$-torus $T$ for which $y(T)$ is a given element of $H^1(F,W)$. 

\subsection{Construction of $T$} We start with a $1$-cocycle of $\Gamma$ in
$W$ and are going to construct $T$ such that $y(T)$ is the class of this
$1$-cocycle. Since the Galois action on $W$ is trivial, the $1$-cocycle is
simply a homomorphism $\Gamma \to W$.

For any integer $l$ not divisible by the characteristic of $k$ we let $F_l$
denote the subfield of
$F_s$ obtained by adjoining an $l$-th root $\epsilon^{1/l}$ of $\epsilon$ to
$F$, and we fix a primitive $l$-th root $\zeta_l$ of $1$ in $k$. We do this
in such a way that $\zeta^m_{lm}=\zeta_l$ for all positive integers $l,m$ not
divisible by the characteristic. Recall that $\Gal(F_l/F)$ is cyclic of order
$l$, generated by the automorphism $\tau_l$ sending $\epsilon^{1/l}$ to
$\zeta_l\epsilon^{1/l}$. 

Since $|W|$ is not divisible by the characteristic, the given homomorphism
$\Gamma \to W$ comes from a homomorphism $\varphi:\Gal(F_{|W|}/F) \to W$, and
we denote by $w$ the element
$\varphi(\tau_{|W|}) \in W$. Choose an element $\dot{w} \in N_G(A)(k)$ of
finite order such that
$\dot{w} \mapsto w$. We denote by $l$ the order of $\dot{w}$. Note that $l$
is not divisible by the characteristic. Thus there is a homomorphism
$\Gal(F_l/F) \to N_G(A)(k)$ sending $\tau_l$ to
$\dot{w}$, and this homomorphism lifts the given homomorphism
$\Gal(F_{|W|}/F) \to W$.   

The element $\dot{w} \in G(k)$ is semisimple, so
there exists a maximal torus $A'$ in $G$ such that $\dot{w}\in A'(k)$. We
have lifted our $1$-cocycle in $W$ to a $1$-cocycle in $N_G(A)$ that 
takes values in $A'$. In order to produce $T$ we need to write this
$1$-cocycle as the coboundary of an element in $G(F_s)$.

But our $1$-cocycle happens to take values in the split torus $A'$, so by
Hilbert's Theorem 90, it must be possible to construct an element $b \in
A'(F_l)$ such that $b^{-1}\tau_l(b)=\dot{w}$. To obtain an explicit element
$b$, we note that since $\dot{w}$ is an $l$-torsion element in $A'(k)$, there
exists $\mu \in X_*(A')$ such that $\mu(\zeta_l^{-1})=\dot{w}$. Now put
$b=\mu(\epsilon^{1/l})^{-1} \in A'(F_l)$; an easy calculation shows that
$b^{-1}\tau_l(b)=\dot{w}$, as desired.

Define a maximal torus $T$ by $T:=bAb^{-1}$. Since $b^{-1}\tau_l(b)=\dot{w}$,
it follows that $T$ is defined over $F$ and that $y(T)$ is the class in
$H^1(F,W)$  we started with. This completes the explicit construction of $T$.

\subsection{Comparison of two filtrations on $\Cartan$} \label{sub4.4} 
The
appropriate context here is that of apartments for tamely
ramified tori in reductive groups over non-archimedean local fields (see
\cite{adler98,yu01} for example).  

We
continue with the maximal torus we just constructed. We use the maximal torus
$A'$ and the element $x:=\mu/l \in X_*(A') \otimes_{\mathbf Z}\mathbf R$ to
obtain an $\mathbf R$-filtration
$\bla_{x,s}$ ($s \in \mathbf R$). It is one of the Moy-Prasad filtrations
considered before, but for $A'$ rather than $A$.

Let $\Cartan $ denote the Lie algebra of~$T$, a Cartan subalgebra of $\bla$.
We want to understand the filtration on $\Cartan$ induced by the filtration
$\bla_{x,s}$ on $\bla$. To do so we first extend the filtration
$\bla_{x,s}$ on $\bla$ to a filtration on $\la(F_l):=\la \otimes_k F_l$
by putting
\[
\la(F_l)_x(s):=\prod_{(\lambda,m)} \la_\lambda \cdot (\epsilon^{1/l})^m
\] where the index set for the product is the set of $(\lambda,m) \in X^*(A')
\oplus \mathbf Z$ such that $\lambda(x)+(m/l) = s$, or equivalently
$\langle \lambda,\mu \rangle+m = sl$, and then putting
\[
\la(F_l)_{x,s}:=\prod_{s' \ge s} \la(F_l)_x(s').
\]

Similarly we  extend the filtration
$\bla_{0,s}$ on $\bla$ to a filtration on $\la(F_l)$
by  putting
\[
\la(F_l)_0(s):= \la \cdot (\epsilon^{1/l})^m
\] 
for $s$  of the form $m/l$ for  $m \in \mathbf Z$ and setting it equal
to $0$ otherwise,  and then putting
\[
\la(F_l)_{0,s}:=\prod_{s' \ge s} \la(F_l)_0(s').
\] 

The group $k^\times$ acts on $F_l$ by automorphisms over $k$. We denote by
$\tau_z \in \Aut(F_l/k)$ the automorphism obtained from $z \in k^\times$; it
sends $\epsilon^{1/l}$ to $z\epsilon^{1/l}$. This action induces an action $z
\mapsto \mathrm{id}_{\la} \otimes \tau_z$ of $k^\times$ on
$\la(F_l)$, and we also denote this action by $\tau_z$. Using $\mu$, we get a
new action $\tau'_z$ of $k^\times$ on $\la(F_l)$, for which $z \in k^\times$
acts by $\Ad(\mu(z)) \circ \tau_z$. For the action $\tau'_z$ an element $z
\in k^\times$ acts on $\la_\lambda \cdot (\epsilon^{1/l})^m$ by
$z^{\langle\lambda,\mu \rangle+m}$, which shows that $\la(F_l)_{x,s}$ is the
product of the weight spaces for the $\tau'_z$ action for weights $j \in
\mathbf Z$ such that $j \ge sl$. Similarly $\la(F_l)_{0,s}$ is the
product of the weight spaces for the $\tau_z$ action for weights $j \in
\mathbf Z$ such that $j \ge sl$.

 An easy calculation shows that
\[
\tau'_z=\Ad(b) \circ \tau_z \circ \Ad(b)^{-1}.
\] Therefore $\Ad(b)^{-1} \la(F_l)_{x,s}=\la(F_l)_{0,s}$. Moreover
$b^{-1}Tb=A$.  This shows that when we use $\Ad(b)^{-1}$ to identify
$\Lie(T_{F_l})$ with $\Lie(A_{F_l})$, the subspaces $\la(F_l)_{x,s} \cap
\Lie(T_{F_l})$ and $\la(F_l)_x(s) \cap \Lie(T_{F_l})$ go over to 
$\la(F_l)_{0,s} \cap \Lie(A_{F_l})$ and $\la(F_l)_0(s) \cap \Lie(A_{F_l})$
respectively. Note that 
$\la(F_l)_0(s) \cap \Lie(A_{F_l})$ is $\Lie(A) \cdot(\epsilon^{1/l})^m$ if
$s=m/l$ for some $m \in
\mathbf Z$ and is $0$ otherwise.

\subsection{Equivalued elements} We continue with $T$ as above. Let
$\val:F_l^\times \to \mathbf Q$ be the valuation on $F_l$, normalized so that
$\val(\epsilon^{1/l})=1/l$. Let $u \in \Cartan$ and suppose that $u$ is
regular, so that $\alpha(u) \in F_l^\times$ for every root $\alpha$ of $T$
(over $\bar F$). Suppose further that there exists $s \in \mathbf Q$ such that
$\val \alpha(u) =s$ for every root $\alpha$ of $T$ (over  $\bar F$) and $\val
\lambda (u) \ge s$ for every $\lambda \in X^*(T)$, in which case we say that
$u$ is  \emph{equivalued} with valuation
$s$.  (Here we are using the convention that $\val 0=+\infty$.)

 This hypothesis on $u$ can also be rephrased as saying that
$u'=\Ad(b)^{-1}u \in
\prod_{m \ge ls}
\Lie(A) \cdot (\epsilon^{1/l})^m$ and that the projection of $u'$ on $\Lie(A)
\cdot (\epsilon^{1/l})^{ls}$ is of the form $u'' \cdot (\epsilon^{1/l})^{ls}$
with $u''$ regular in
$\Lie(A)$. It follows from our comparison of filtrations in \ref{sub4.4} that
$u \in \bla_{x,s}$ and that the image of $u$ under
\[
\bla_{x,s} \to \bla_x(s)=\la(x,\mathbf s) \hookrightarrow \la,
\] namely $u''$, is regular semisimple in $\la$, hence is $G$-good in $\la$
(see \ref{good.vect}, noting that if $W$ is non-trivial, its order is
divisible by $2$, so that characteristic $2$ has been excluded). 
 
A minor nuisance is that $x$ lies in $X_*(A')\otimes_{\mathbf Z}\mathbf R$
rather than $\apt$. We can correct this by replacing $T,u$ by their
conjugates under an element of $G(k)$ that takes $A'$ to $A$. This replaces
$x$ by an element of $\apt$, which we will still call $x$. For this new $x$
and for any $t \in \mathbf R$ with $t \le s$ we have shown that the three
conditions of
\ref{sub3.2} are satisfied. We conclude that for any $y \in \apt$ and $t \le
s$ the intersection of the affine Springer fiber $\mathcal F_y(t,u)$ with any
$\bg_x$-orbit on $\mathcal P_y$ is an iterated affine space bundle over a
Hessenberg variety for $\bg_x(0)$. (The Hessenberg variety depends on the
orbit.)

\subsection{Weakly Coxeter tori} For convenience we now assume that the
maximal $F$-split torus $A_T$ in~$T$ is contained in $A$, so that the
centralizer
$M$ of $A_T$ is a Levi subgroup of $G$ (over $k$) containing $A$. We further
assume that $T$ is weakly Coxeter, so that $w$ is a Coxeter element in $W_M$,
the absolute Weyl group of~$A$ in~$M$. Our lifting $\dot w$ is a regular
semisimple element in~$M$ by a result of Kostant \cite{kostant59}. We are free
to modify $\dot w$ by any element in the center of~$M$; doing so we may
choose $\dot w$ so that it is regular in~$G$, not just in~$M$. Looking back
at how $x$ was defined, we see that $x$ lies in the interior of an alcove.
Thus $\bg_x$ is an Iwahori subgroup, and the Hessenberg varieties we get are
subvarieties of the flag manifold of the torus $\bg_x(0)$ and hence have
either $1$ or $0$ elements. Thus, for weakly Coxeter $T$ 
 our Hessenberg
paving is a paving by affine spaces.

\section{Regular semisimple implies good  when
 $p \ne 2$} In this section we fix a prime number $p$ and assume that our
algebraically closed ground field $k$ has  characteristic $p$.

\subsection{Review of regularity in $\mathfrak g$}
Recall from \cite[Prop. 11.8]{borel69} that $u \in \mathfrak g$ is semisimple
if and only if it is contained in the Lie algebra of some maximal torus $T$
of~$G$, in which case it is  clear that the centralizer $\mathfrak g_u$ 
of~$u$ in~$\mathfrak g$ contains the Lie algebra of that maximal torus. A
semisimple element $u \in \mathfrak g$ is said to be \emph{regular} if its
centralizer is as small as possible, or, in other words, is equal to the Lie
algebra of some maximal torus. We denote by $\mathfrak g_{\rs}$ the set of
regular semisimple elements in~$\mathfrak g$. 

\subsection{Non-emptiness of $\mathfrak g_{\rs}$}

 Consider elements $u$ in the Lie algebra $\mathfrak t$ of some fixed maximal
torus~$T$ in~$G$. Then $u$ is regular if and only $\alpha(u) \ne 0$ for every
root $\alpha$ of~$T$ in~$G$. By conjugacy of maximal tori,  $\mathfrak
g_{\rs}$ is non-empty if and only if the set of regular semisimple elements in
$\mathfrak t$ is non-empty, and this happens if and only if no root of $T$
vanishes identically on~$\mathfrak t$. It even suffices to look only at
simple roots, since any root is in the Weyl group orbit of a simple one. 
(A root vanishes identically on $\mathfrak t$ if and only if it is
divisible  by $p$ in~$X^*(T)$.)
 
If $G$ is adjoint no simple root is divisible by
$p$, and therefore $\mathfrak g_{\rs}$ is non-empty.  On the other hand
$\mathfrak g_{\rs}$ is in fact empty for
$G=Sp_{2n}$ and $p=2$. Since $\langle \alpha^\vee, \alpha \rangle =2$,  a
root is never divisible by a prime other than $2$; therefore 
$\mathfrak g_{\rs}$ is non-empty if $p \ne2$.

\subsection{Chevalley's restriction theorem} Chevalley's theorem states that
so long as 
$\mathfrak g_{\rs}$ is non-empty, the restriction map from $G$-invariant
polynomials on $\mathfrak g$ to Weyl group invariant polynomials
on~$\mathfrak t$ is an isomorphism. We refer to Springer and
Steinberg \cite[Section
$3.17$]{springer-steinberg70} for a proof that works in finite
characteristic.  (Springer and Steinberg assume that $G$ is adjoint, but it
is routine to adapt their proof to the more general case in which it is
assumed only that 
$\mathfrak g_{\rs}$ is non-empty.)

\subsection{Parallel results for $\mathfrak g^*$}
Consider the dual $\mathfrak g^*$ of $\mathfrak g$, a representation of~$G$.
Since $\mathfrak t$ is canonically a direct summand of~$\mathfrak g$ (with
complement the direct sum of the root subspaces for~$T$), we may identify the
dual $\mathfrak t^*$ of $\mathfrak t$ with a  direct
summand of~$\mathfrak g^*$. For $\lambda \in \mathfrak t^*$ it is easy to see
that the centralizer of $\lambda$ in $\mathfrak g$ (for the coadjoint action
of $G$ on $\mathfrak g^*$) is the direct sum of $\mathfrak t$ and the root
subspaces $\mathfrak g_\alpha$ for roots $\alpha$ of~$T$ such that
$\lambda(H_\alpha)=0$. Here $H_{\alpha}:=\alpha^\vee(1)$, in other words, the
image of the coroot $\alpha^\vee$ in~$\mathfrak t$.
An element $\lambda \in \mathfrak t^*$ is said to be \emph{regular} if its
centralizer is as small as possible, or, in other words, is equal to
$\mathfrak t$.  We denote by $\mathfrak g^*_{\rs}$ the set of
 elements in~$\mathfrak g^*$ whose coadjoint orbit contains a regular element
of $\mathfrak t^*$. 

Considerations parallel to those above show that 
$\mathfrak g^*_{\rs}$ is non-empty if $p\ne 2$ or if $G$ is semisimple simply
connected. On the other hand $\mathfrak g^*_{\rs}$ is empty for $SO(2n+1)$
and $p=2$.

Finally, it is routine to transpose the proof of Springer and Steinberg
from $\mathfrak g$ to $\mathfrak g^*$, where
it yields the statement that if
$\mathfrak g^*_{\rs}$ is non-empty, then the restriction map from
$G$-invariant polynomials on $\mathfrak g^*$ to Weyl group invariant
polynomials on~$\mathfrak t^*$ is an isomorphism.

\subsection{Goodness of elements in $\mathfrak g_{\rs}$ and  $\mathfrak
g_{\rs}^*$} \label{good.proof}
Assume that both  $\mathfrak g_{\rs}$ and  $\mathfrak
g_{\rs}^*$ are non-empty. Then Chevalley's restriction theorem holds for both
$\mathfrak g$ and $\mathfrak g^*$, and it is therefore clear that there are
no non-zero unstable vectors (in the sense of \ref{unstab.vect}) 
 in
$\mathfrak t$ and
$\mathfrak t^*$. It follows easily that elements in $\mathfrak g_{\rs}$ and
$\mathfrak g^*_{\rs}$ are good vectors (in the sense of \ref{good.vect}).

\bibliographystyle{amsalpha}

\providecommand{\bysame}{\leavevmode\hbox to3em{\hrulefill}\thinspace}

\end{document}